\documentclass[12pt]{article}

\usepackage{pstricks}
\usepackage{amsmath,amsthm}
\usepackage{amssymb}
\usepackage{color}
\usepackage{xspace}
\usepackage[colorlinks=true,
linkcolor=green,
filecolor=brown,
citecolor=green]{hyperref}
\definecolor{gray}{rgb}{.221,.221,.221}

\def\red{\textcolor{red} }
\def\v{\vert}

\def\lr{LRmax\xspace}
\def\rr{revert-reciprocal\xspace}

\begin{document}
\newtheorem{theorem}{Theorem}
\newtheorem{defn}[theorem]{Definition}
\newtheorem{lemma}[theorem]{Lemma}
\newtheorem{prop}[theorem]{Proposition}
\newtheorem{cor}[theorem]{Corollary}
\begin{center}

{\Large
Lagrange Inversion Counts $3\overline{5}241$-Avoiding Permutations                          \\ 
}

\vspace{10mm}
DAVID CALLAN  \\
Department of Statistics  \\
\vspace*{-2mm}
University of Wisconsin-Madison  \\
\vspace*{-2mm}
1300 University Ave  \\
\vspace*{-2mm}
Madison, WI \ 53706-1532  \\
{\bf callan@stat.wisc.edu}  \\
\vspace{2mm}

April 19, 2011
\end{center}

\begin{abstract}
In a previous paper, we showed that $3\overline{5}241$-avoiding permutations are counted 
by the unique sequence that starts with a 1 and shifts left under 
the self-composition transform. The proof uses a complicated bijection. 
Here we give a much simpler proof based on Lagrange inversion.
\end{abstract}

\section{Introduction}  
A permutation avoids the barred pattern $3\overline{5}241$ when the (not necessarily consecutive)
pattern 3241 occurs only as part of a 35241 pattern.
The composition of two sequences $(a_{n})_{n\ge1}$ and $(b_{n})_{n\ge1}$ 
is defined by composition of their
(ordinary) generating functions. A sequence is unital if its first term is 1.
There is a unique unital sequence $(b_{n})_{n\ge1}$ whose composition with itself 
is equal to its left shift, $(b_{2},b_{3},\ldots)$, the so-called 
left-shift eigensequence for self-composition. This sequence 
begins $(b_{n})_{n\ge1}=(1,1,2,6,23,104,531,\ldots)$.
In \cite{eseq06}, we showed that the counting sequence for 
$3\overline{5}241$-avoiding permutations is this sequence indexed 
from 0. More precisely, with 
$S_{n}(\tau)$ denoting the set of permutations of $[n]$ that avoid 
the pattern $\tau$, we established
\begin{theorem}\label{mainresult}
    The sequence $(a_{n})_{n\ge1}$ defined by $a_{1}=1$ and, for 
    $n\ge1,\ a_{n+1}=\v\,S_{n}(3\overline{5}241)\,\v$ is the left-shift 
eigensequence for self-composition.
\end{theorem}
The proof used a rather complicated bijection. Our objective here is to simplify the 
proof using Lagrange inversion. We start with the characterization of 
$3\overline{5}241$-avoiding permutations given in \cite{eseq06}, and recall some notation. 
Every permutation $\pi$ on $[n]$, considered as a list (or word), has the form 
$\pi=m_{1}L_{1}m_{2}L_{2} \ldots m_{r}L_{r}$ where 
$m_{1}<m_{2}<\ldots<m_{r}=n$ are the left-to-right maxima  of $\pi$ (\lr for 
short) and the $L_{i}$ are the intervening words. 
We call $m_{1}L_{1}\,/\,m_{2}L_{2}\,/\,\ldots \,/\, 
m_{r}L_{r}$ the \lr decomposition of $\pi$.
\begin{theorem} \emph{\cite[Theorem 1]{eseq06}}\label{characterize}
   A permutation $\pi$ on $[n]$ is $3\overline{5}241$-avoiding if and only 
   if its \lr decomposition satisfies
    
    $($\mbox{\,i}\,$)$ $L_{1}<L_{2}<\ldots<L_{r}$ in the sense that $u\in 
    L_{i},\,v\in L_{j}$ with $i<j $ implies $u<v$, and 
    
    $($ii\,$)$ each $L_{i}$ is $3\overline{5}241$-avoiding. \qed   
\end{theorem}

\vspace*{6mm}

\section{The Revert-Reciprocal Transform}  \label{revrec}
Let us define a transform on sequences, $(a_{n})_{n\ge1} \rightarrow (b_{n})_{n\ge1}$, 
using generating functions, by
\[
A(x) \rightarrow B(x):=
\big( \textrm{\raisebox{0.4mm}{$\frac{x}{1+A(x)}$}} 
\big)^{\textrm{{\footnotesize $\langle -1 \rangle$}}},
\]
where ${}^{\langle -1 \rangle}$ denotes compositional inverse.
We'll call it the \emph{\rr}transform.
The Lagrange inversion formula \cite[Thm. 5.4.2, p.\,38]{ec2} is just the ticket to find an explicit 
form for the entries of the transformed sequence---we find that $b_{1}=1$ and for $n\ge 1$,
\[
    b_{n+1}=\frac{1}{n+1} \sum_{1^{r_{1}}\ldots\: n^{r_{n}}\vdash n}
\binom{n+1}{n+1-{\scriptstyle \sum}_{i=1}^{n}r_{i},\:r_{1},\: \ldots\:,\: r_{n}}
a_{1}^{r_{1}}a_{2}^{r_{2}}\ldots a_{n}^{r_{n}},
\]
where the sum is over all partitions $1^{r_{1}}\ldots n^{r_{n}}$ of 
$n$, written in frequency-count form. Thus the unique fixed point 
$(a_{n})_{n\ge 1}$ for 
the \rr transform is defined by $a_{1}=1$ and for $n\ge 1$,
\begin{equation} \label{fixptseq}
a_{n+1}=\frac{1}{n+1} \sum_{1^{r_{1}}\ldots \: n^{r_{n}}\vdash n}
\binom{n+1}{n+1-{\scriptstyle \sum}_{i=1}^{n}r_{i},\:r_{1},\: \ldots\:,\: r_{n}}
a_{1}^{r_{1}}a_{2}^{r_{2}}\ldots a_{n}^{r_{n}}.
\end{equation}
\begin{prop}\label{FP=eseq}
    The fixed point for \rr coincides with the left-shift 
eigensequence for self-composition.
\end{prop}
Proof.\quad With $B(x):=
\big( \textrm{\raisebox{0.4mm}{$\frac{x}{1+A(x)}$}} 
\big)^{\textrm{{\scriptsize $\langle -1 \rangle$}}}$, the 
\rr transform of $A(x)$, we have
\begin{eqnarray*}
    B(x) \hspace*{6mm} & = & A(x)  \hspace*{5mm} \Leftrightarrow  
    \\[2mm]
    \frac{A(x)}{1+A(A(x))} & = & \hspace*{2mm} x \hspace*{9mm} \Leftrightarrow  \\
    A(A(x))\hspace*{3mm} & = & 
    \frac{A(x)-x}{\textrm{\raisebox{1mm}{$x$}}}, 
\end{eqnarray*}
which is the defining relation for the left-shift 
eigensequence for self-composition. \qed

\vspace*{6mm}

\section{Permutations and Trees}  
We consider a permutation as a word of distinct letters from the 
alphabet of positive integers. A \emph{cycle} is a permutation in which the 
largest entry occurs first. A \emph{standard} permutation is one on an 
initial segment of the positive integers. To \emph{standardize} a permutation 
means to replace its smallest entry by 1, second smallest by 2, and so on.
A \emph{cycle-labeled} ordered tree is an ordered tree in which, for each non-leaf 
vertex $v$, the child edges of $v$ are labeled, left to right, with the entries of a 
standard cycle.
\begin{theorem}\label{permAsTree}
    Permutations on $[n]$ satisfying condition $($\mbox{i}\,$)$ of 
    Theorem \ref{characterize} are in bijective correspondence with 
    cycle-labeled ordered trees on $n$ edges.
\end{theorem}
Proof.\quad In the \lr decomposition, $m_{1}L_{1}\: / \:m_{2}L_{2}\: / \:\ldots \: / \: 
m_{r}L_{r}$, of a permutation $\pi$ on $[n]$, 
set $a_{i}=$ length of $m_{i}L_{i}$ and $b_{i}=m_{i}-m_{i-1}$ for 
$1\le i\le r$, with $m_{0}:=0$. Now create $a_{r}$ edges from the root 
labeled with the entries of $m_{r}L_{r}$ from left to right. Then, at the 
$b_{r}$-th leaf, place 
$a_{r-1}$ edges labeled with $m_{r-1}L_{r-1}$. Proceed to place $a_{r-2}$ 
edges labeled with $m_{r-2}L_{r-2}$ at the $b_{r-1}$-th leaf in the 
current tree, and so on. The result for $\pi=3\:1\:2\:/\:5\:4\:/\ 11\ 7\ 6\ 8\:/\:
12\:/\ 14\ 13\ 10\ 9$ is shown on the left in Figure 1 below.

\begin{center}
\begin{pspicture}(-2,-2.1)(4,5)
\psset{unit=.7cm}
\newrgbcolor{gren}{0 1 0}

\psline (-4,2)(0,0)(-2,2)(-2,4)(-4.5,6)
\psline(-3,6)(-2,4)(0.5,6)
\psline(-1,6)(-2,4)
\psline(1,2)(0,0)(3,2)(2,4)
\psline(3,2)(4,4)(2.5,6) 
\psline(4,6)(4,4)(5.5,6)

\psdots(-4,2)(0,0)(-2,2)(-2,4)(-4.5,6)(-3,6)(-2,4)(0.5,6)(-1,6)(1,2)(0,0)(3,2)(2,4)(4,4)(2.5,6)(4,6)(4,4)(5.5,6) 

\rput(1,-2){\textrm{{\small $(a_{i})_{i=1}^{r}=(3,2,4,1,4),\ \ (b_{i})_{i=1}^{r}=(3,2,6,1,2)$}}}

\rput(-4.2,5.3){\textrm{{\footnotesize $11$}}}
\rput(-2.4,5.3){\textrm{{\footnotesize $7$}}}
\rput(-1.6,5.3){\textrm{{\footnotesize $6$}}}
\rput(0,5.3){\textrm{{\footnotesize $8$}}}
\rput(-2.4,3){\textrm{{\footnotesize $12$}}}
\rput(-3.2,1.2){\textrm{{\footnotesize $14$}}}

\rput(-0.7,1.2){\textrm{{\footnotesize $13$}}}
\rput(1,1.2){\textrm{{\footnotesize $10$}}}
\rput(2.3,1.2){\textrm{{\footnotesize $9$}}}
\rput(2.1,3.2){\textrm{{\footnotesize $5$}}}
\rput(3.9,3.2){\textrm{{\footnotesize $4$}}}
\rput(3.8,5.3){\textrm{{\footnotesize $1$}}}

\rput(2.7,5.3){\textrm{{\footnotesize $3$}}}
\rput(5.3,5.3){\textrm{{\footnotesize $2$}}}

\rput(0,-1){\textrm{{\small edge-labeled tree}}}

\end{pspicture}
\begin{pspicture}(-4.5,-2.1)(3,5)
\psset{unit=.7cm}

\psline (-4,2)(0,0)(-2,2)(-2,4)(-4.5,6)
\psline(-3,6)(-2,4)(0.5,6)
\psline(-1,6)(-2,4)
\psline(1,2)(0,0)(3,2)(2,4)
\psline(3,2)(4,4)(2.5,6) 
\psline(4,6)(4,4)(5.5,6)

\psdots(-4,2)(0,0)(-2,2)(-2,4)(-4.5,6)(-3,6)(-2,4)(0.5,6)(-1,6)(1,2)(0,0)(3,2)(2,4)(4,4)(2.5,6)(4,6)(4,4)(5.5,6) 

\rput(-4.2,5.3){\textrm{{\red{\footnotesize $4$}}}}
\rput(-2.4,5.3){\textrm{{\footnotesize $2$}}}
\rput(-1.6,5.3){\textrm{{\footnotesize $1$}}}
\rput(0,5.3){\textrm{{\footnotesize $3$}}}
\rput(-2.3,3){\textrm{{\red{\footnotesize $1$}}}}
\rput(-3.2,1.2){\textrm{{\red{\footnotesize $4$}}}}

\rput(-0.7,1.2){\textrm{{\footnotesize $3$}}}
\rput(1,1.2){\textrm{{\footnotesize $2$}}}
\rput(2.3,1.2){\textrm{{\footnotesize $1$}}}
\rput(2.1,3.2){\textrm{{\red{\footnotesize $2$}}}}
\rput(3.9,3.2){\textrm{{\footnotesize $1$}}}

\rput(3.8,5.3){\textrm{{\footnotesize $1$}}}
\rput(2.7,5.3){\textrm{{\red{\footnotesize $3$}}}}
\rput(5.3,5.3){\textrm{{\footnotesize $2$}}}

\rput(0,-1.4){\textrm{{\small cycle-labeled tree}}}


\end{pspicture}
\end{center}

\vspace*{-5mm}

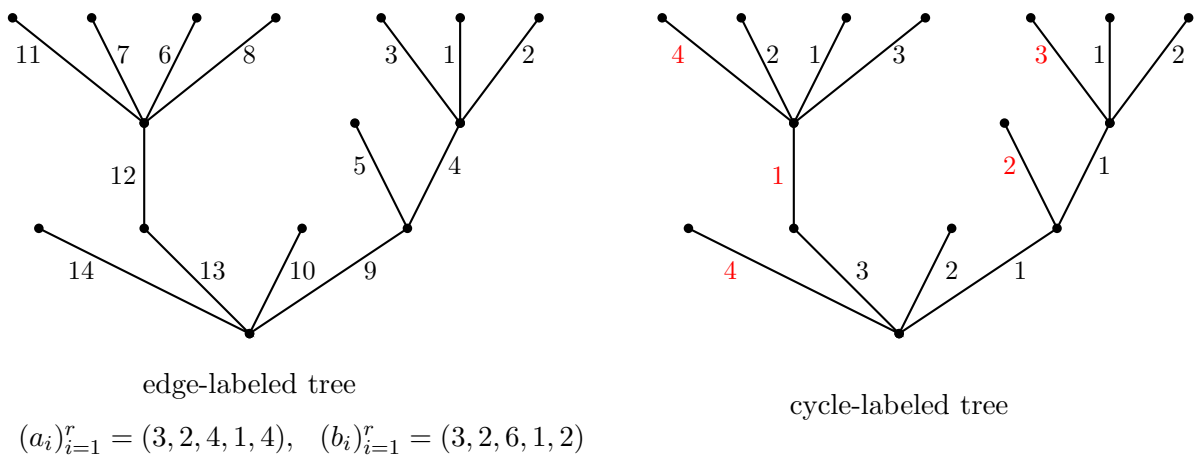
\centerline{Figure 1}
\vspace*{2mm}
\noindent Next, for each non-leaf vertex, erase the label on its leftmost child 
edge (a \lr entry), standardize the permutation labeling the remaining   
child edges, and relabel the leftmost child edge with the outdegree 
of the vertex, thereby obtaining a standard cycle as the labeling on the child 
edges of the vertex. 
The result is shown on the right in Figure 1 above with leftmost child 
edge labels in red. 
 
The \lr entries are encoded in the cycle-labeled tree by the locations of the 
non-leaf vertices, and so can be recovered. The leftmost child edge 
labels on non-leaf vertices taken in preorder (aka clockwise 
walkaround order) give the lengths of the $m_{i}L_{i}$ segments in 
reverse order. Since we know 
the \lr entries, the support sets of the $L_{i}$'s can then be determined. 
Finally, the permutations labeling the child edges of the non-leaf 
vertices (ignoring the label on the leftmost edge) determine the 
actual words $L_{i}$ from their support sets. \qed

\section{The Coup de Gr\^{a}ce} 
The number of ordered trees on $n$ edges (hence, $n+1$ vertices) with 
outdegree sequence $(r_{i})_{i=0}^{n}$---$r_{i}$ is the number of 
vertices with $i$ children---is \cite[Thm 5.3.10,\:p.\,34]{ec2}
\[
\frac{1}{n+1}\binom{n+1}{r_{0},r_{1},\ldots,r_{n}}.
\]
Necessarily, $\sum_{i=0}^{n}r_{i} = n+1$, the number of vertices, and 
$\sum_{i=0}^{n}ir_{i} = n$, the number of edges.
Under the correspondence of Theorem \ref{permAsTree}, Theorem 
\ref{characterize} (ii) says that $3\overline{5}241$-avoiding permutations 
on $[n]$ correspond to cycle-labeled trees with $n$ edges in which each cycle (or, 
equivalently, each cycle without its first entry) is 
$3\overline{5}241$-avoiding. 

Consequently, by counting cycle-labeled trees in which the label 
lists associated with the vertices are $3\overline{5}241$-avoiding permutations, 
we find, letting $a_{n}$ denote the 
number of $3\overline{5}241$-avoiding permutations on $[n-1]$, that
\[
a_{n+1}=\frac{1}{n+1} \sum_{1^{r_{1}}\ldots\: n^{r_{n}}\vdash n}
\binom{n+1}{n+1-{\scriptstyle \sum}_{i=1}^{n}r_{i},\:r_{1},\: \ldots\:,\: r_{n}}
a_{1}^{r_{1}}a_{2}^{r_{2}}\ldots a_{n}^{r_{n}}.
\]
But this is precisely the relation (\ref{fixptseq}) that defines the 
fixed point sequence for \rr and so, since this fixed 
point coincides with the left-shift eigensequence for self-composition 
by Prop. \ref{FP=eseq}, Theorem \ref{mainresult} is established.

\section{Concluding Remarks}
The \rr transform as defined in Section \ref{revrec} is convenient for 
our purposes. But the transform always starts with a 1. If we 
redefine the \rr transform to delete this 1, 
\[
A(x) \rightarrow B(x):=
\frac{\left( \textrm{\raisebox{0.8mm}{$\frac{\textrm{\raisebox{0.8mm}{\normalsize 
$x$}}}{\textrm{\raisebox{-0.4mm}{\normalsize $1+A(x)$}}}$}} 
\right)^{\textrm{{\footnotesize $\langle -1 \rangle$}}}-1}{x},
\]
then it becomes an invertible
transform from sequences $(a_{n})_{n\ge1}$ to $(b_{n})_{n\ge1}$ and, 
with this modified definition of \rr, Prop. \ref{FP=eseq} has the form

\emph{The (unique) left-shift eigensequences for the transforms \rr and 
self-composition coincide.
}

\noindent The two transforms, however, are quite different: the only sequence on 
which they agree is this left-shift eigensequence. Self-composition 
is not quite invertible, but it is invertible on the set of real 
sequences $\{(a_{n})_{n\ge 1}\,:\ a_{1}>0\}$.

\end{document}